\documentclass[a4paper,10pt]{article}
\usepackage[utf8]{inputenc}

\usepackage{graphicx}
\usepackage{times}
\usepackage{amssymb}
\usepackage{amsmath}
\usepackage{cite}

\usepackage[font={small, sf, bf}, labelfont=bf,labelsep=period, textfont=md]{caption}
\usepackage[a4paper, left=25mm, right=20mm, top=20mm, bottom=25mm, nohead]{geometry}

\makeatletter
\renewcommand{\@biblabel}[1]{\quad#1.}
\makeatother

\title{The probability of drawing intersections: extending the hypergeometric distribution}
\author{Alex T. Kalinka\footnote{alex.t.kalinka@gmail.com}}
\date{\small{Institute for Population Genetics, Vetmeduni, Veterin\"{a}rplatz 1, Vienna, Austria.}}

\begin{document}

\maketitle 

\begin{abstract}

The wide availability of biological data at the genome-scale and across multiple variables has resulted in statistical questions regarding the enrichment or depletion of the number of discrete objects (e.g. genes) identified in individual experiments. Here, I consider the problem of inferring enrichment or depletion when drawing independently, and without replacement, from two or more separate urns in which the same $n$ distinct categories of objects exist. The statistic of interest is the size of the intersection of object categories. I derive a probability mass function describing the distribution of intersection sizes when sampling from $N$ urns and show that this distribution follows the classic hypergeometric distribution when $N=2$. I apply the theory to the intersection of genes belonging to a set of traits in three different vertebrate species illustrating that the use of $P$-values from one-tailed enrichment tests enables accurate clustering of related traits, yet this is not possible when relying on intersection sizes alone. In addition, intersection distributions provide a means to test for co-localization of objects in images when using discretized data, allowing co-localization tests in more than two channels. Finally, I show how to extend the problem to variable numbers of objects belonging to each category, and discuss how to make further progress in this direction. The distribution functions are implemented and freely available in the R package `hint'.

\end{abstract}

\section{Introduction}

\noindent Biological data can be highly multivariate and recent advances in high-throughput technologies has increased the ease with which such data can be acquired \cite{suetal2004}. When handling data which is composed of observations made on sets of discrete objects, such as genes or proteins, it may often be necessary to ask whether the number of objects identified in a particular treatment is greater or less than expected by chance. If objects can be classified into two categories and sampling is without replacement from a single urn, then the hypergeometric distribution suffices to describe the distribution of the number of objects belonging to a particular category.

\begin{figure}[!ht]
\begin{center}
\advance\leftskip-0.1cm
\includegraphics[width=4.5in]{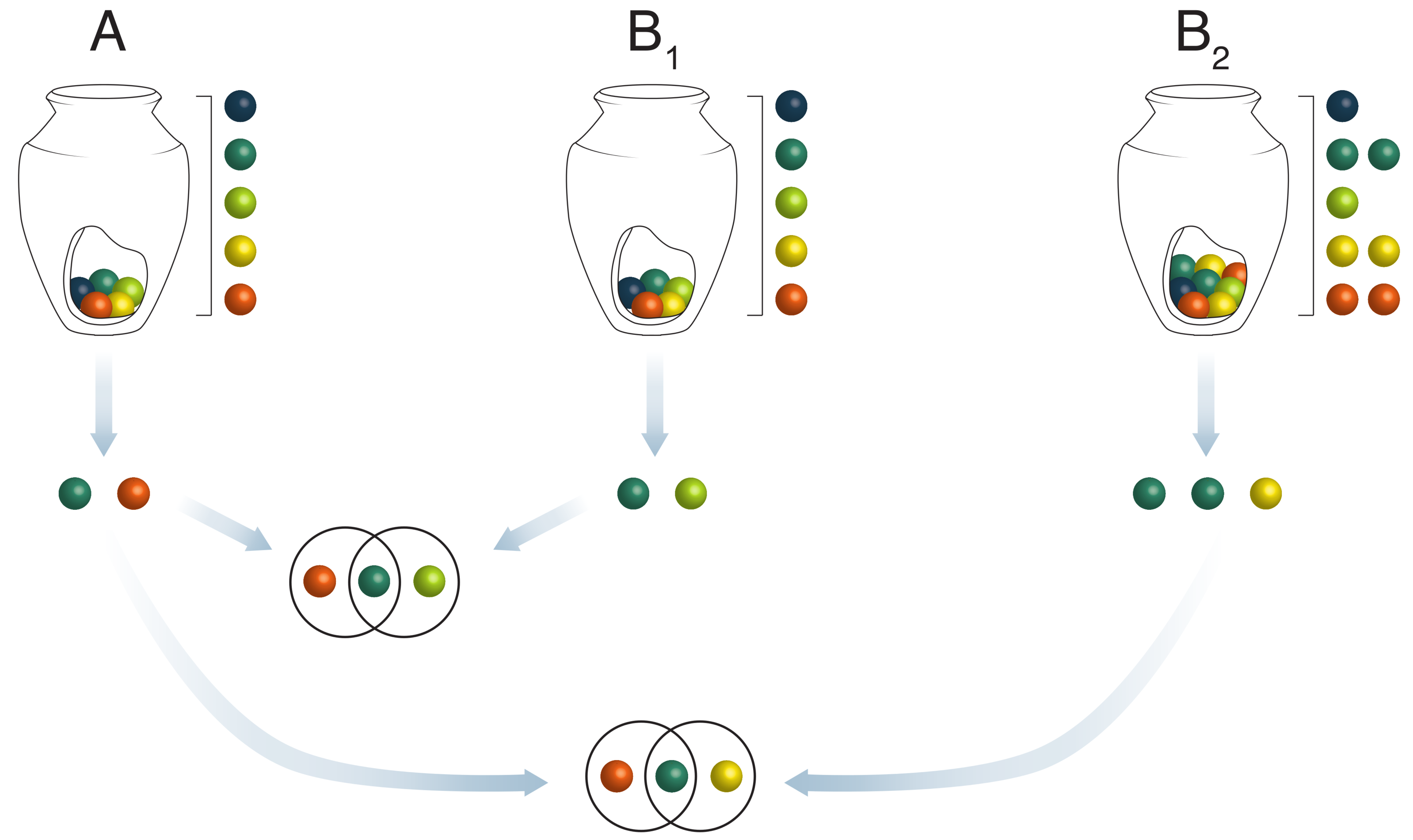}
\end{center}
\caption{
A schematic illustrating the drawing of intersections from urns containing balls belonging to 5 different categories (depicted using different colours). In urns \textbf{A} and \textbf{B$\mathbf{_1}$} there is exactly 1 ball in each of the categories, whereas in urn \textbf{B$\mathbf{_2}$} 3 of the categories contain duplicate members. Although both duplicates of one category are drawn from this urn, the intersection size remains 1.
}
\label{urns_fig}
\end{figure}

However, we may instead be interested in the intersection of object categories when sampling independently from several urns (Figure \ref{urns_fig}). Such a scenario might arise if we know, for example, that a set of genes in one species shares a particular characteristic (e.g. expression in a particular organ) in common with a set of genes in another species; then we might ask whether there is a significant enrichment of homologous genes in the intersection of both gene sets when categories are defined by homology \cite{heynetal2014}.

The hypergeometric distribution \cite{johnsonetal2005} describes the probability of $k$ successes in $n$ draws without replacement from a single population of size $N$ in which reside $D$ possible successes, and is given by

\begin{equation} \label{hypdist}
\mathbb{P}(X = k) = \frac{ \binom{D}{k} \binom{N-D}{n-k} }{\binom{N}{n} }.
\end{equation}

\noindent The distribution has been broadly applied to tests of significance for categorical data in which objects can be classified in two different ways \cite{pearson1899, fisher1922, gonin1936}.

Imagine instead that we have two separate urns each containing objects that belong to one of $n$ distinct categories. If we draw $a$ objects from the first urn and $b$ objects from the second, what is the probability of finding an intersection of size $v$ in the categories drawn from both urns? Thus, in contrast to the hypergeometric distribution, this problem involves more than two categories of objects and independent sampling from two, or more, separate urns (Figure \ref{urns_fig}).

In what follows, I derive probability mass functions describing the distribution of intersection sizes when sampling from two or more urns. I show that in the case of two urns, the distribution is hypergeometric. This result illustrates that the hypergeometric distribution can be used to describe sampling from two urns in addition to the classic single urn interpretation. I then use the distributions to infer the relationships between biological traits among three vertebrate species and demonstrate that the use of enrichment tests uncovers the true relationships between the traits. The remainder is devoted to extending the approach to allow for variable numbers of objects in each of the $n$ categories.

\section{Symmetrical, singleton case}

\subsection{Two urns}

\noindent First, I consider the simplest scenario in which there are two urns containing exactly one member in each of the $n$ categories, i.e. a symmetrical, singleton case (corresponding to sampling from urns \textbf{A} and \textbf{B$\mathbf{_1}$} in Figure \ref{urns_fig}). We sample $a \leq n$ and $b \leq n$ from each urn respectively and wish to know the probability of drawing intersections of size $v$ where

\[ \max(a+b-n,0) \leq v \leq \min(a,b). \]

To count the number of ways of picking an intersection of size $v$, it is useful to note that we can count the number of ways of picking a single, specific combination of intersecting categories (e.g. categories \{1,2,3\} for $v=3$) by counting the number of non-intersecting categories that can be drawn to produce this particular intersection combination. Hence, for the first urn there are $\binom{n-v}{a-v}$ ways to draw one particular combination of intersecting categories. This leaves $\binom{n-a}{b-v}$ ways of drawing from the second urn to give an intersection size of $v$ for a single, specific combination of categories; the upper index of $(n-a)$ ensures that we do not count intersections of size larger than $v$. The total number of ways to pick intersections of size $v$ must then be summed over all $\binom{n}{v}$ category combinations:

\[ C_{v} = \sum^{\binom{n}{v}} \binom{n-v}{a-v}\binom{n-a}{b-v} = \binom{n}{v}\binom{n-v}{a-v}\binom{n-a}{b-v}. \]

\noindent The probability of picking an intersection of size $v$ is then $C_{v}$ divided by the total number of ways of picking $a$ and $b$ from $n$:

\begin{equation} \label{symmsingprecursor}
\mathbb{P}(X = v) = \frac{ \binom{n}{v}\binom{n-v}{a-v}\binom{n-a}{b-v} }{ \binom{n}{a}\binom{n}{b} }.
\end{equation}

\noindent Applying a trinomial revision \cite{grahametal1994} to the first two binomials in the numerator, the expression can be reduced to

\begin{equation} \label{symmsingdist}
\mathbb{P}(X = v) = \frac{ \binom{n}{a}\binom{a}{v}\binom{n-a}{b-v} }{ \binom{n}{a}\binom{n}{b} } = \frac{ \binom{a}{v}\binom{n-a}{b-v} }{\binom{n}{b}},
\end{equation}

\noindent which is the hypergeometric distribution given in Equation \ref{hypdist}, and is symmetrical in terms $a$ and $b$. The symmetry of the problem, in which both urns contain exactly 1 member in each of the $n$ categories, enables this simplification. This derivation illustrates that the hypergeometric distribution can be used to describe sampling from two urns as well as the classic single urn interpretation. Simulations show that the distribution is exact (Suppl. Figure \ref{simhyp_fig}).

\subsection{$N$ urns}

\noindent When sampling is from $N > 2$ urns, we need to account for intersections between fewer than $N$ urns (among the non-intersecting categories) since they will contribute to our statistic of interest, $v$, which measures intersections across all of the $N$ urns. However, for each cross-urn intersection between less than $N$ urns, it is sufficient to account for intersections between $N-1$ of the urns since the problem is fully specified by $N-1$ urns. Hence, for intersections between $k$ urns, there will be $\binom{N-1}{k}$ cross-urn intersections that must be accounted for. To arrive at the total number, we must sum over all cross-urn intersections smaller than $N$:

\[ \sum^{N-1}_{k=2} \binom{N-1}{k} = 2^{N-1}-N .\]

\noindent Thus, when there are three urns, we must account for intersections between the two urns that belong to the $N-1$ urns. For example, with three urns, \{$A,B,C$\}, we would need to consider intersections between urns $A$ and $B$. The maximum intersection size in this case is $\alpha = \min(a-v,b-v)$.

When summing over all possible pair-wise intersections, each draw that is shared by both $A$ and $B$ in the non-intersecting categories (not intersecting across all $N$ urns) will be drawn from the $a-v$ non-intersecting categories drawn from $A$, and to avoid double counting these shared categories, we must subtract them from the $b-v$ drawn from $B$:

\[ \sum_{i=0}^{\alpha} \binom{a-v}{i} \binom{n-a}{b-v-i} . \]

\noindent To ensure that these pair-wise intersections are not counted in the $c-v$ items drawn from $C$ (since they are intersections across $A$ and $B$ only), they must be subtracted from
the $n-v$ items that can be drawn from $C$, giving:

\[ \sum_{i=0}^{\alpha} \binom{a-v}{i} \binom{n-a}{b-v-i} \binom{n-v-i}{c-v}. \]

\noindent The probability of drawing an intersection of size $v$ across all three urns is then

\begin{equation} \label{hypint3}
 \mathbb{P}(X = v|N=3) = \frac{ \binom{a}{v} \sum_{i} \binom{a-v}{i} \binom{n-a}{b-v-i} \binom{n-v-i}{c-v} }{ \binom{n}{b} \binom{n}{c} } .
\end{equation}

\noindent Simulations confirm that the distribution is exact (Suppl. Figure \ref{multi_3}). Although a closed form for the above expression is not readily apparent, the sum in the numerator can be re-arranged, following the procedure outlined by Roy \cite{roy1987} and Hirschhorn \cite{hirschhorn2002}, into the following expression

\[ \binom{n-a}{b-v} \binom{n-v}{c-v} {}_3F_2 \left( \begin{matrix}c-n,v-a,v-b \\ v-n, 1+n+v-a-b\end{matrix} \hspace{0.125cm} \middle\vert 1 \right) \]

\noindent where ${}_{3}F_{2}(.)$ denotes the generalized hypergeometric function ($(.)_i$ denotes the rising factorial):

\[ \sum^{\infty}_{i=0} \frac{(c-n)_i(v-a)_i(v-b)_i}{i!(v-n)_i(1+n+v-a-b)_i}, \]

\noindent which, since it is neither balanced nor well-poised \cite{roy1987}, does not appear to permit a closed form \cite{milgram2010}, although the series will terminate when $i > \alpha$. Nonetheless, this implies the following interesting identity

\begin{equation} \label{hyper_identity}
 \sum_{v\geq0} \binom{a}{v} \binom{n-a}{b-v} \binom{n-v}{c-v} {}_3F_2 \left( \begin{matrix}c-n,v-a,v-b \\ v-n, 1+n+v-a-b\end{matrix} \hspace{0.125cm} \middle\vert 1 \right) = \binom{n}{b} \binom{n}{c}.
\end{equation}
 
Extending to the case of four urns will require that we account for three pair-wise urn intersections ($A:B, A:C, B:C$) together with one three-way intersection ($A:B:C$) ($2^{3}-4=4$). Following the logic above, we can derive the following expression for the number of ways of picking $v$ intersections across four urns:
 
\[ \binom{a}{v} \sum_{i,j,k,l} \binom{a-v}{i} \binom{a-v-i}{j-l} \binom{b-v-i}{k-l} \binom{i}{l} \binom{n-a}{b-v-i} \binom{n-a-b+v+i}{c-v-j-k+l} \binom{n-v-l}{d-v} \]

\noindent where $i,j,k,l$ represent sums over intersections in $A:B$, $A:C$, $B:C$, and $A:B:C$ respectively and $d$ is the number sampled from the fourth urn, $D$. This can be simplified by applying Vandermonde convolutions to the sums in $j,k$ to give the following distribution (simulations shown in Suppl. Figure \ref{multi_4}):

\begin{equation}
 \mathbb{P}(X=v|N=4) = \frac{ \binom{a}{v} \sum_{i,l} \binom{a-v}{i}\binom{i}{l}\binom{n-a}{b-v-i}\binom{n-v-i}{c-v-l}\binom{n-v-l}{d-v} }{ \binom{n}{b}\binom{n}{c}\binom{n}{d} }.
\end{equation}

\noindent By observing that the the four nested sums have reduced to two sums over the first cross-urn intersections in each intersection group ($A:B$ and $A:B:C$; ordering urns from $A$ to $D$), we can infer that this will also be the case for $N$ urns. Hence, for 5 urns, we will have 3 sums over $A:B$, $A:B:C$, and $A:B:C:D$ intersections, meaning that there will be $N-2$ nested sums for $N$ urns. Thus, for the general case, we have the following probability mass function:

\begin{equation} \label{general_hypergeom}
 \mathbb{P}(X=v|N) = \left. { \binom{a_{1}}{v} \sum_{i_{j} \in S} \prod_{j=0}^{|S|} \binom{i_{j}}{i_{j+1}} \binom{n-v-i_{j}}{a_{j+2}-v-i_{j+1}} } \middle/ {\prod_{k=2}^{N} \binom{n}{a_{k}} } \right.
\end{equation}

\noindent where the sum is a nested sum over the set, $S$, of $N-2$ cross-urn intersections, $a_{k}$ is the number drawn from the $k$'th urn, and $i_{j}$ takes the following values when the index $j$ is beyond the indices of members of $S (j = 1...N-2)$:

\[  i_{j} = 
 \begin{cases}
  a_{1}-v, & \text{if } j = 0 \\
  0, & \text{if } j > N-2  \\
 \end{cases} \]

\noindent Simulations for $N = 5,6$ confirm that the distribution is exact (Suppl. Figures \ref{multi_5},\ref{multi_6}). Substituting $N=2$ into the above expression, we can recover the classic hypergeometric distribution, showing that it is a special case of this more general distribution. It is worth noting that this distribution implies the identity of the sum over $v$ of the numerator of Equation \ref{general_hypergeom} with the denominator, which is a generalisation of the identity given in Equation \ref{hyper_identity}.

It is possible to deduce the first two moments of this distribution since we know that the distribution is hypergeometric for $N=2$ urns. The expectation for the hypergeometric is $ab/n$, and hence we can infer that the expectation for $N$ urns is

\begin{equation}
 \mathbb{E}(X) = \frac{\prod_{k=1}^{N} a_{k} }{n^{N-1}}.
\end{equation}

\noindent This can be confirmed numerically \cite{kalinka2013}, and can also be derived using a binomial approximation (see below). Knowing the expectation can help us to infer the variance using the relationship $\mathbb{V}\mathrm{ar}(X) = \mathbb{E}(X^2)-\mathbb{E}(X)^2$. For the hypergeometric case ($N=2$), we have

\[ \frac{ab(n+ab-a-b)}{n(n-1)} - \left(\frac{ab}{n}\right)^2. \]

\noindent The expression contained in the brackets of the numerator of the first fraction can be written as $(n-1)+(1-a)(1-b)$, and, hence, we can infer the general expression for $N$ urns:

\begin{equation}
 \mathbb{V}\mathrm{ar}(X) = \frac{ (\prod_{k=1}^{N} a_{k}) ((n-1)^{N-1} + (-1)^{N} \prod_{k=1}^{N} (1-a_{k}))  }{n^{N-1}(n-1)^{N-1}} - \left(\frac{\prod_{k=1}^{N} a_{k} }{n^{N-1}}\right)^{2}
\end{equation}

\noindent which can again be confirmed numerically \cite{kalinka2013}. Knowing the expectation and variance for the general case enables the use of Gaussian approximations for calculating probabilities when the number of categories and the sample sizes are all large \cite{nicholson1956}. We can also see that the following is true

\[ \forall a_{k} < n : \lim_{N \to \infty} \mathbb{E}(X) = 0 \]

\noindent and since the lower bound of the distribution is 0, the limit must also be true for the variance. Hence, we have the intuitive result that as the number of urns grows large, the probability of picking an intersection across all of them tends to 0.

\subsection{Approximation for large $n$}

\noindent The binomial distribution is a good approximation for the hypergeometric distribution when the total population is large ($N$ in Equation \ref{hypdist}) and the sample drawn from this population is small ($n$ in Equation \ref{hypdist}). Therefore, it is reasonable to suppose that a similar approximation will apply to intersections across $N$ urns. I start by extracting Equation \ref{hypint3}, and assume that $b$ is small relative to $n$, and that $a$ and $c$ grow large as $n$ grows large. Multiplying top and bottom by $(b-v)!, (n-v)!$ and $(n-v-i)!$, the expression can be re-written as

\[ \binom{b}{v} \frac{c!(n-v)!}{(c-v)!n!} \sum_{i=0}^{b-v} \binom{b-v}{i} \frac{a!(n-v-i)!}{(a-v-i)!n!} \frac{(n-a)!(n-b)!}{(n-a-(b-v-i))!(n-b+(b-v-i))!} \frac{(n-c)!(n-v-i)!}{(n-c-i)!(n-v)!}. \]

\noindent If we further note that 

\[ \frac{c!(n-v)!}{(c-v)!n!} = \frac{c^{\underline v}}{n^{\underline v}} = \prod_{k=1}^{v} \frac{c-v+k}{n-v+k} \]

\noindent and that

\[ \lim_{n \to \infty} \prod_{k=1}^{v} \frac{c-v+k}{n-v+k} = \prod_{k=1}^{v} \lim_{n \to \infty} \frac{c-v+k}{n-v+k} = \left(\frac{c}{n}\right)^{v} = p_{c}^{v} \]

\noindent then the expression reduces to

\[ \binom{b}{v} p_{c}^{v} \sum_{i=0}^{b-v} \binom{b-v}{i} p_{a}^{v+i} (1-p_{a})^{b-v-i} (1-p_{c})^{i}.\]

\noindent The sum can be evaluated using the binomial theorem if we take $p_{a}^{v}$ outside

\[ \binom{b}{v} p_{c}^{v} p_{a}^{v} \sum_{i=0}^{b-v} \binom{b-v}{i} (1-p_{a})^{b-v-i} (p_{a}-p_{a}p_{c})^{i} = \binom{b}{v} p_{c}^{v} p_{a}^{v} (p_{a}-p_{a}p_{c}+1-p_{a})^{b-v} = \binom{b}{v} (p_{a}p_{c})^{v} (1-p_{a}p_{c})^{b-v}. \]

\noindent Hence, the approximation for $N$ urns can be readily deduced as

\begin{equation} \label{binomint}
 \mathbb{P}(X = v|N) = \binom{b}{v} \left(\prod_{i=1}^{N-1} p_{i}\right)^{v} \left(1 - \prod_{i=1}^{N-1} p_{i}\right)^{b-v}
\end{equation}

\noindent which will hold when $n$ is large and the samples from $N-1$ urns are larger than the sample from one of the urns (only one urn has a small sample). The distribution is a variant of the binomial and could be fairly described as a binomial intersection distribution. The expectation and variance are easily derived as

\[ \mathbb{E}(X) = b \prod_{i=1}^{N-1} p_{i}  \hspace{1cm}\text{and}\hspace{1cm} \mathbb{V}\mathrm{ar}(X) = b \prod_{i=1}^{N-1} p_{i} \left(1 - \prod_{i=1}^{N-1} p_{i}\right). \]

\noindent This expectation is also the expectation of the true distribution, but the variance is greater (as is the case for the binomial and hypergeometric distributions). From these expressions, it can be seen that in the limit of large $N$ both the expectation and the variance tend to zero:

\[ \lim_{N \to \infty} b \prod_{i=1}^{N-1} p_{i} = 0  \hspace{1cm}\text{and}\hspace{1cm} \lim_{N \to \infty} b \prod_{i=1}^{N-1} p_{i} \left(1 - \prod_{i=1}^{N-1} p_{i}\right) = 0 .\]

\noindent These limits will also hold for the exact distribution when the sample sizes for each urn are not equal to $n$. Numerical comparisons demonstrate that Equation \ref{binomint} is a good approximation for the true distribution (Suppl. Figure \ref{binom_approx}).

\section{Trait relationships across three species}

\begin{figure}[!ht]
\begin{center}
\advance\leftskip-0.25cm
\includegraphics[width=6.5in]{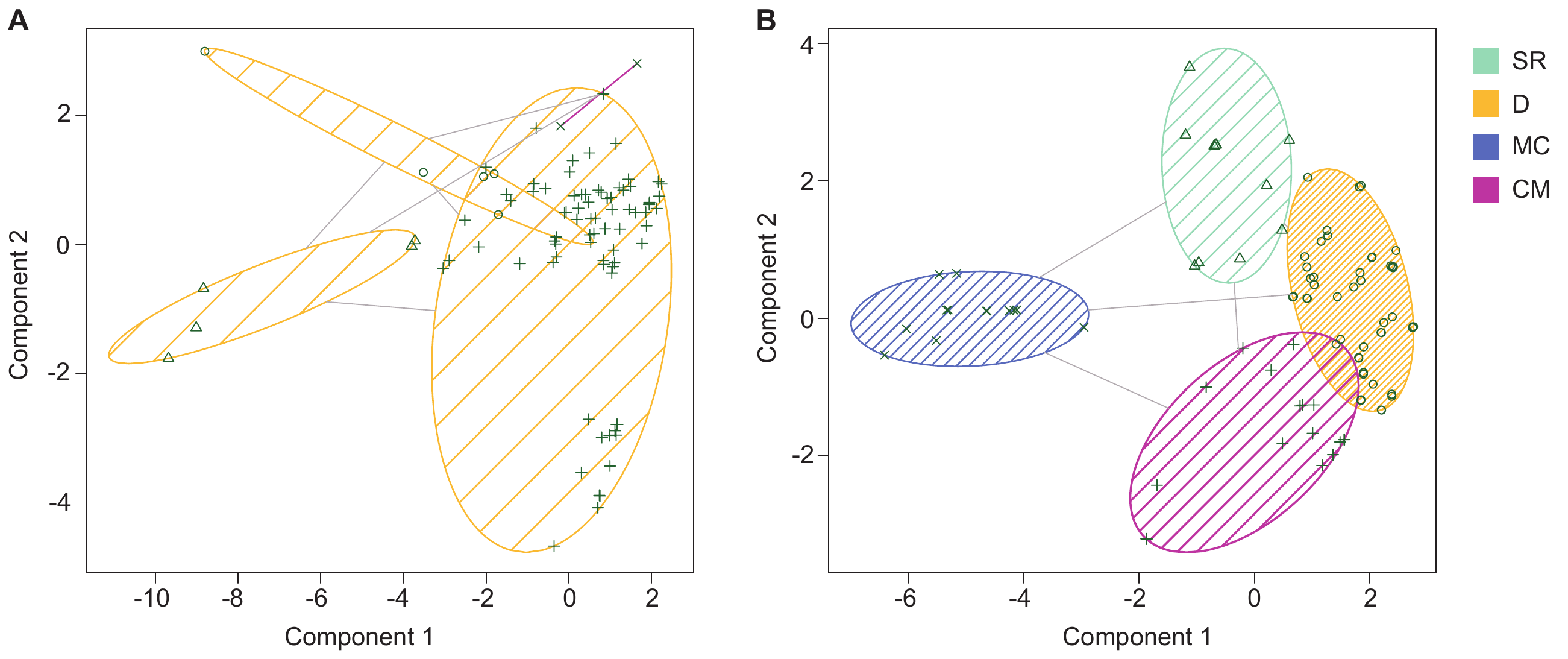}
\end{center}
\caption{
$K$-means clustering of biological traits across three vertebrate species (human, mouse, and zebrafish). In \textbf{A} traits were clustered using intersection sizes and in \textbf{B} traits were clustered using $P$-values from one-tailed enrichment tests based on Equation \ref{hypint3}. SR: sexual reproduction, D: development, MC: metabolism core, CM: cabohydrate metabolism.
}
\label{clusters}
\end{figure}

\noindent To demonstrate the utility of these distributions, I used them to test for enrichment of orthologous genes across three vertebrate species (human, mouse, and zebrafish). For each species, I downloaded genes from Ensembl Biomart \cite{kinsellaetal2011} that were annotated to a set of traits based on Gene Ontology (GO) biological function terms \cite{ashburneretal2000}. Traits were chosen so that they fell into four well-defined categories: development, sexual reproduction, carbohydrate metabolism, and core metabolism. The number of one-to-one orthologs that were shared across all three species was determined pair-wise by traits (e.g. brain development vs oogenesis) and across all three species combinations (e.g. human trait 1 vs mouse trait 2 vs zebrafish trait 1, etc).

The resulting matrix was then clustered using the $K$-means clustering algorithm \cite{steinhaus1957, macqueen1967} implemented in R \cite{R2012}. Clustering was conducted separately using either the intersection size of one-to-one orthologs across trait and species comparisons, or $P$-values of one-tailed enrichment tests based on Equation \ref{hypint3} and implemented in the R package `hint' \cite{kalinka2013}. The results illustrate that when clustering by intersection size alone, traits cannot be distinguished into separate clusters, and genes belonging to gene-rich traits (developmental traits) fall into three of the four clusters (Figure \ref{clusters}A). In contrast, when clustering according to enrichment tests, the four trait groups can be clearly distinguished from each other (Figure \ref{clusters}B). Intersection tests effectively control for the number of genes shared across species and the number of genes present in the trait as a whole, thereby identifying traits that are highly related even if 
they have relatively small numbers of genes.

Although I have tested for enrichments of one-to-one orthologs across species, it is also possible to test for enrichments within a single species; for example, we could test for enrichment of genes expressed in different organs belonging to a single species. Many other variants are also possible. In addition, intersection distributions provide a means to test for significant co-localization of objects in digital images in which two or more distinct fluourescent labels have been imaged \cite{mandersetal1993}. Fluouresent intensity data will need to be discretized by applying a cutoff below which the presence of a label is considered within background levels. One advantage is that Equation \ref{general_hypergeom} allows the overlap of any number of labels (or channels) to be tested.

\section{Asymmetrical cases}

\subsection{Duplicates in one of two urns}

\noindent If we allow duplicates in $q \leq n$ of the categories in the second urn (each category can contain 1 or 2 balls but not 0), then the problem becomes asymmetrical (corresponding to \textbf{A} and \textbf{B$\mathbf{_2}$} in Figure \ref{urns_fig}). The presence of duplicates in the second urn will reduce the overall chance of drawing an intersection of a certain size because duplicates that are both sampled from the urn can at most contribute an intersection of size 1 (see Figure \ref{urns_fig}). Thus, it is reasonable to conjecture that the expectation for this distribution will always be less than for the equivalent symmetrical, singleton case:

\begin{equation}\label{conjecture}
\mathbb{E}(X | q \neq 0) < \frac{a b}{n}.
\end{equation}

\noindent Furthermore, we can reasonably suppose that the expression describing this distribution will be a variant of the hypergeometric since we have made only a small modification to the basic problem, and added a single parameter, $q$. I begin, therefore, by modifying Equation \ref{symmsingprecursor} to account for the effects of including $q$ duplicates.\\

\noindent There are three main differences affecting the drawing of both intersecting and non-intersecting categories:

\begin{enumerate}
\item{When a category is sampled from the first urn (among the $a-v$ non-intersecting categories) for which there is a duplicate pair in the equivalent category in the second urn, this reduces the number of ways we can pick non-intersecting categories from the second urn to ensure an intersection size of $v$. The number of such draws is indicated by the index $m$.}
\item{If a category with a duplicate pair is picked in the $v$ intersecting items, this does not reduce the number of ways of picking non-intersecting items from the second urn (since drawing the duplicate member will also produce an intersection of size $v$), but it removes a duplicate category from the $m$ that could be picked in the non-intersecting set. The number of such draws is indicated by the index $l$.}
\item{Picking $l$ duplicate categories in the $v$ intersecting items increases the number of ways that these items can be drawn, but for each duplicate that is picked in $v$, one less duplicate is available for the non-intersecting set to be drawn from the second urn. The number of such draws is indicated by the index $j$.}
\end{enumerate}

\noindent To calculate the probability, we must sum over all of the ways of combining the above events such that they produce intersection sizes of $v$, which must satisfy

\[ \max\left(a+b-n-\min\left(\left\lfloor \frac{b}{2} \right\rfloor,q\right),0\right) \leq v \leq \min(a,b) \]

\noindent where the lower bound is determined by the maximum number of duplicates that can be picked from the second urn. I will move from left to right across the numerators of Equation \ref{symmsingprecursor} and describe how each binomial term must be modified. The number of ways of drawing $v$ intersecting categories, $\binom{n}{v}$, must incorporate consideration for the $l$ duplicate categories listed in point 2 above, leading to:

\[ \sum_{l}\sum^{l}_{j}\binom{n-q}{v-l}\binom{q}{l}\binom{l}{j}, \]

\noindent which counts the total number of ways of picking $v$ with and without $l$ duplicates. The number of ways of picking non-intersecting categories from the first urn for a single category combination, $\binom{n-v}{a-v}$, must be modified to account for the $m$ duplicate categories listed in point 1 above:

\[ \sum_{m}\sum_{l}\binom{n-v-q+l}{a-v-m}\binom{q-l}{m}, \]

\noindent which counts the number of ways of picking $a-v$ non-intersecting categories given that we have sampled both $m$ and $l$ duplicates. Finally, the number of ways of picking non-intersecting items from the second urn, $\binom{n-a}{b-v}$, must be modified to account for a reduction in duplicates that can be drawn to ensure an intersection size of $v$:

\[ \sum_{m}\sum^{l}_{j}\binom{n+q-a-m-j}{b-v}, \]

\noindent which counts the number of ways of picking $b-v$ non-intersecting categories from the second urn given that we have picked $m$ non-intersecting duplicate equivalents and $j$ intersecting duplicate equivalents from the first urn.

Summing over all the possible combinations of these events then gives us the total number of ways of picking an intersection of size $v$ in the duplicate case (underbraces indicate equivalent expressions in the symmetrical singleton case in Equation \ref{symmsingprecursor}):

\[
 C^{d}_{v} = \sum_{m=0}^{\beta} \sum_{l=0}^{\gamma} \sum_{j=0}^{l} \underbrace{ \binom{n-q}{v-l} \binom{q}{l} \binom{l}{j} }_{\binom{n}{v}} \underbrace{ \binom{q-l}{m} \binom{n-v-q+l}{a-v-m} }_{\binom{n-v}{a-v}} \underbrace{ \binom{n+q-a-m-j}{b-v} }_{\binom{n-a}{b-v}},
\]

\noindent where

\[ \beta = \min(a-v,q) \hspace{1cm}\text{and}\hspace{1cm} \gamma = \min(v,q-m). \]

\noindent The probability is then $C^{d}_{v}$ divided by the total number of ways of picking $a$ and $b$ from both urns:

\begin{equation}
 \mathbb{P}(X=v) = \frac{ \sum_{m,l,j} \binom{n-q}{v-l} \binom{q}{l} \binom{l}{j} \binom{q-l}{m} \binom{n-v-q+l}{a-v-m}  \binom{n+q-a-m-j}{b-v} }{ \binom{n}{a}\binom{n+q}{b} }
\end{equation}

\noindent A closed-form expression for the above equation is not forthcoming. Simulations show that the distribution is exact (Suppl. Figure \ref{simhyp_dup_fig}). This derivation illustrates that a small change in the details of the urn model (allowing duplicates in one urn) greatly complicates the form of the probability mass function.

\section{The distribution of intersection distances}

\begin{figure}[!ht]
\begin{center}
\advance\leftskip-1.75cm
\includegraphics[width=6in]{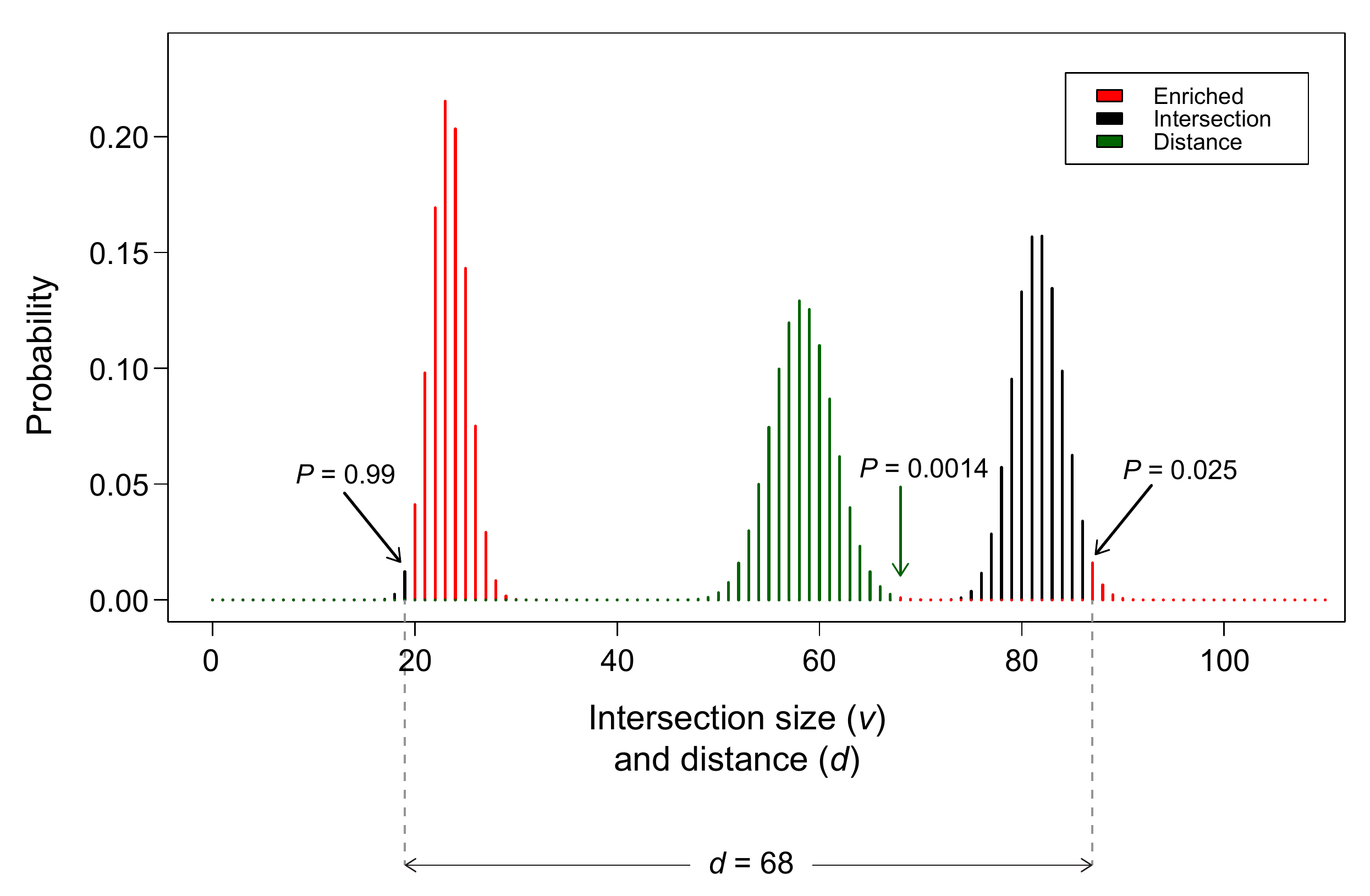}
\end{center}
\caption{
Two intersection distributions (in black and red) and their distance distribution (in dark green and red). One-tailed tests for greater intersection sizes than expected by chance have been applied to both intersection distributions at 19 (left) and 87 (right). This gives a distance of 68, which is greater than expected by chance ($P = 0.0014$) even though significance at the 5\% level was marginal for only one of the intersection distributions. Parameter values for the two intersection distributions are, from left to right: $n=60,a=40,b=35; n=155,a=110,b=115$.
}
\label{dist_dist}
\end{figure}

\noindent When drawing intersections from two different distributions (with different parameters, or, for example, with a singleton case and a duplicate case) it might be of interest to ask whether the absolute distance between their intersection sizes is what would be expected by chance. Testing for significant differences between intersection sizes is likely to be of interest when we want to know if they are behaving differently, i.e. are the intersection sizes that we observe falling into opposite tails more than would be expected by chance (Figure \ref{dist_dist})? In biological terms, the question would be whether we have strong enrichment and depletion in a pair of observations.

To calculate the probability of finding an intersection distance of size $d$, we need to sum over all the ways to produce $d$ when pairing our two distributions:

\begin{equation}
\mathbb{P}(X=d) = \sum_{\{v_{1},v_{2}\}_{i} \in D_{d}}^{|D_{d}|} \mathbb{P}(v_{1_i}|n_{1},a_{1},...)\cdot\mathbb{P}(v_{2_i}|n_{2},a_{2},...)
\end{equation}

\noindent where $D_{d}$ is the set of pairs of intersection sizes, $\{v_{1},v_{2}\}$, with absolute differences of size $d$. If $R$ and $S$ are the sets of all possible intersection sizes for both distributions, then

\[  |D_{d}| = 
 \begin{cases}
  |R \cap S|, & \text{if } d = 0 \\
  \min(|R|, |S|-d) + \min(|S|, |R|-d), & \text{if } d>0 \\
 \end{cases} \]

\noindent where $\min(|R|, |S|-d),\min(|S|, |R|-d)\geq0$ (i.e. negative values are set to zero). This distribution has a relationship to the intersection distributions that is similar to the relationship between the binomial and Bernoulli distributions.

\section{Related distributions: drawing from a single urn}

\noindent It is useful to consider the related, though simpler, distribution associated with drawing $a$ balls from $n$ categories with $q\leq n$ duplicates from a single urn. In this case, we are no longer interested in intersection sizes, but rather in the number of distinct categories, $c$, which are drawn from the single urn. When $q = 0$ then $c = a$ necessarily. Thus, we restrict ourselves to cases where $0 < q \leq n$. The bounds on $c$ are then

\[ a - \min\left(\left\lfloor \frac{a}{2} \right\rfloor,q\right) \leq c \leq a \]

\noindent where the lower bound is determined by the maximum possible number of duplicate pairs subtracted from $a$. For any particular value of $c$, there are always $a-c$ duplicate pairs that must be picked to ensure that there are $c$ distinct categories drawn.

To count the number of ways of drawing $c$ categories, it is useful to first note that there are 3 combinations that need to be counted:

\begin{enumerate}
 \item The number of ways of picking $a-c$ duplicate pairs, i.e. \{1,1\}, \{2,2\}.
 \item The number of ways of picking duplicates not picked as pairs from the $q-a+c$ remaining.
 \item The number of ways of picking non-duplicates.
\end{enumerate}

\noindent Point 1 is simply given by $\binom{q}{a-c}$. For points 2 and 3, we must sum over all the ways of combining duplicates (not picked as pairs) and non-duplicates to give $c$ distinct categories. An important quantity here is the number of non-duplicate pairs (not \{1,1\} or \{2,2\}) in $a$, given by $a - 2(a-c) = 2c-a$. Then

\[ \sum^{q}_{j=0} \binom{q-a+c}{j}\binom{n-a+c-j}{2c-a-j} \]

\noindent gives the combined number for points 2 and 3. The probability of picking $c$ distinct categories is then given by

\begin{equation}
 \mathbb{P}(X=c) = \frac{ \binom{q}{a-c}\sum^{q}_{j=0} \binom{q-a+c}{j}\binom{n-a+c-j}{2c-a-j} }{ \binom{n+q}{a} }.
\end{equation}

\noindent Again, a closed-form expression is not easily derived. However, if we focus on the special case when $q = n$, the above expression can be simplified. Substituting $q=n$ and applying a trinomial revision to the binomials within the sum, the numerator can be reduced to

\[ \binom{n}{a-c}\binom{n-a+c}{2c-a}\sum^{n}_{j=0} \binom{2c-a}{j}, \]

\noindent which in turn simplifies to

\[ \binom{n}{c}\binom{c}{a-c}2^{2c-a}. \]

\noindent From left to right, the three terms count the number of ways of picking $c$ distinct categories from $n$, the number of ways of picking $a-c$ duplicate pairs from $c$, and the number of ways of picking $2c-a$ duplicates \emph{not} picked as duplicate pairs. The probability of drawing $c$ distinct categories when all categories in the urn contain a duplicate is then given by

\begin{equation}\label{hydist}
 \mathbb{P}(X=c) = \frac{ \binom{n}{c}\binom{c}{a-c}2^{2c-a} }{ \binom{2n}{a} }.
\end{equation}

\noindent Here, we find that the numerator and denominator satisfy identity 3.22 appearing in Gould's compendium of combinatorial identities \cite{gould1972}:

\[ \sum^{a}_{c=\left\lfloor \frac{a}{2} \right\rfloor} \binom{n}{c} \binom{c}{a-c}2^{2c} = 2^{a}\binom{2n}{a} \]

\noindent where $2^a$ cancels since we have $2^{2c}2^{-a}$. Using standard approaches \cite{wilf2005}, we can derive the expectation of the distribution. First, we note that

\[ \sum_{c}^{a} \binom{n-1}{c-1} \binom{c}{a-c} 2^{2c} = \frac{2^a}{n} \binom{2n}{a} \mathbb{E}(X). \]

\noindent Working with the LHS, we derive the generating function from which we can derive the RHS:

\begin{equation*}
\begin{aligned}
& \hspace{4.05mm} \sum_{c}^{a} \binom{n-1}{c-1} 2^{2c} x^{c} \sum_{a} \binom{c}{a-c} x^{a-c} \\
&= \sum_{c}^{a} \binom{n-1}{c-1} 2^{2c} x^{c}(1+x)^{c} \\
&= (4x+4x^{2}) \sum_{c}^{a} \binom{n-1}{c-1} (4x+4x^{2})^{c-1}.
\end{aligned}
\end{equation*}

\noindent Hence, we have the following generating function

\begin{equation*}
\begin{aligned}
g(x) &= (4x+4x^{2})(1+4x+4x^{2})^{n-1}\\
&= 4x(1+x)(1+2x)^{2n-2}
\end{aligned}
\end{equation*}

\noindent from which we extract the coefficient of $x^{a}$

\[ [x^{a}]g(x) = \frac{2^a}{n} \binom{2n}{a} a\left(1-\frac{a-1}{4n-2}\right) \]

\noindent and therefore

\begin{equation}
\mathbb{E}(X) = a\left(1-\frac{a-1}{4n-2}\right)
\end{equation}

\noindent which is simply the sample size multiplied by the probability that no duplicates are picked after $a$ draws ($1-\frac{1}{2}\cdot\frac{a-1}{2n-1}$; defined for $a>0$). We can see that when $a=2n$ we will always sample all $n$ of the categories in the urn. When this is not the case, however, the expected number of categories drawn will always be less than $n$, and less than $a$, thereby providing some support for the conjecture given in Equation \ref{conjecture}. 


Further work on single urn distributions will help to shed light on sampling across several urns when there are variable numbers of balls in each category.

\section*{Summary}

\noindent I have discussed a sampling-without-replacement problem that is closely related to the hypergeometric distribution. The main differences are:

\begin{enumerate}
 \item Samples are taken independently from two or more separate urns.
 \item More than 2 categories of objects are allowed.
 \item The statistic of interest is the size of the intersection across the urns.
\end{enumerate}

When there are two urns with exactly one ball in each of the $n$ categories, I have shown that the distribution is hypergeometric. I have also derived a general expression that describes sampling from $N$ urns and showed that the expectation and variance tend to 0 as $N \to \infty$.

With one small modification to the basic two-urn scenario - the addition of $q$ duplicate categories in the second urn - the problem becomes much more complex, and the distribution no longer can be expressed in a closed form, though it is clear that this asymmetrical case is a variant of the hypergeometric. Furthermore, I derive the distribution of absolute distances between the intersection sizes of two separate intersection distributions. This distribution has utility when the question of interest is whether two intersection sizes are behaving significantly differently. Finally, I derived a closed-form expression for the distribution of distinct categories sampled from a single urn containing duplicates in all of its categories. This related distribution may aid in the understanding of intersection distributions and their properties, and ultimately is an attempt to work towards a more general description of this broad class of distributions. More generally, these results highlight that despite the 
extensive study of univariate discrete distributions \cite{johnsonetal2005}, much may remain to be discovered \cite{baker2000,rodriguezavietal2003,murat&szynal2006,karlis&xekalaki2008,satheeshkumar2009}.

\section*{Acknowledgements}

\noindent Thanks are due to Peter Steinbach for assistance with implementing the duplicate case in C++ for large parameter sets, and Iva Kelava for preparing Figures \ref{urns_fig} and \ref{clusters}.

\bibliography{../../refs/evol_genet}

\clearpage

\section*{Supplementary material}

\setcounter{figure}{0}
\renewcommand{\figurename}{Supplementary Figure}

\begin{figure}[!ht]
\begin{center}
\advance\leftskip-0.1cm
\includegraphics[width=4in]{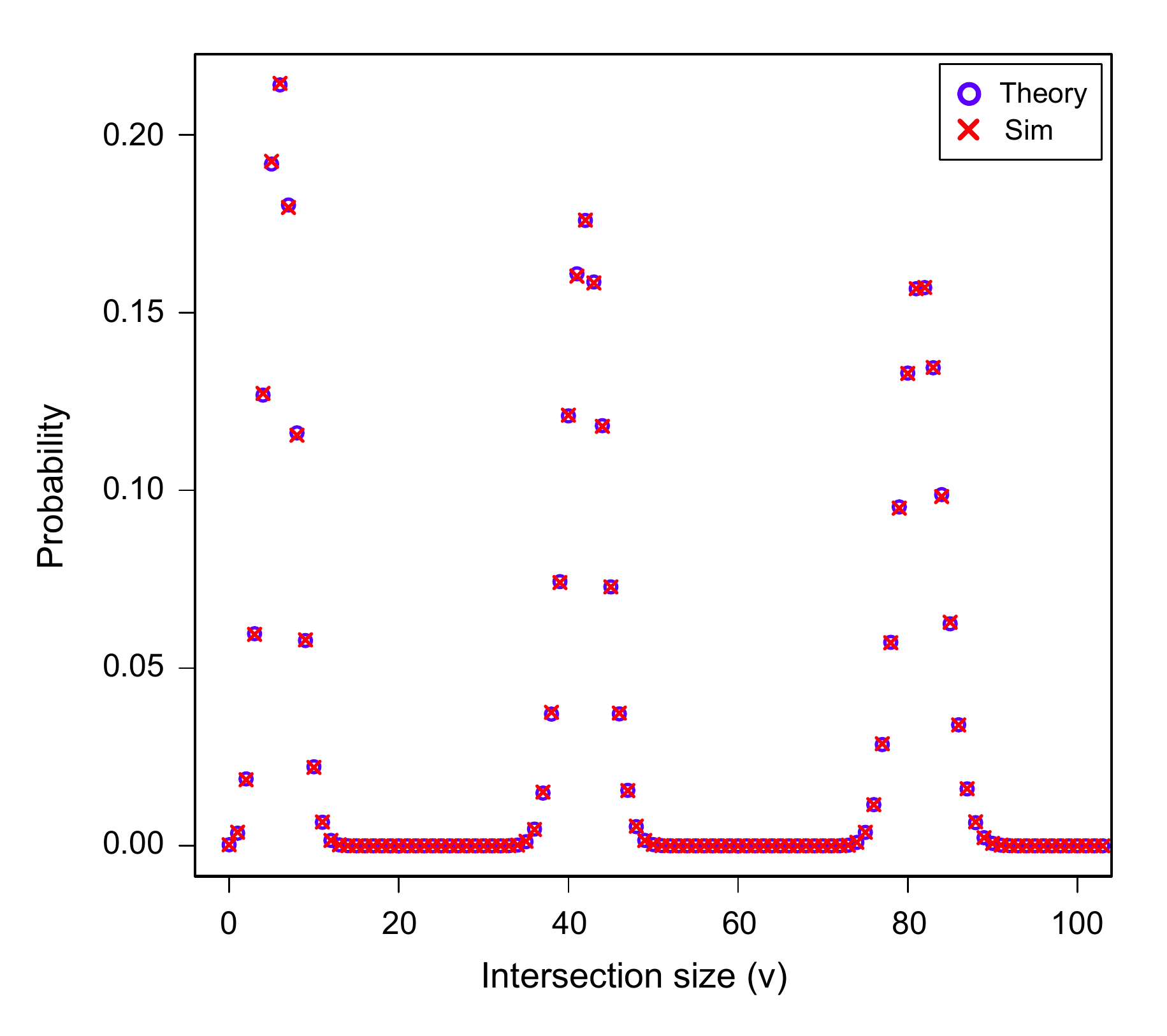}
\end{center}
\caption{
Match between theory and simulation for 3 parameter sets in the symmetrical, singleton case (\textbf{A} and \textbf{B$\mathbf{_1}$} in Figure \ref{urns_fig}). Simulations (Sim) consisted of randomly and independently sampling twice (without replacement) from $n$ distinct categories and recording the size of the intersection each time (repeated 500,000 times for each distribution). From left to right, the parameters were: $n = 100, a = 20, b = 30$; $n=100,a=70,b=60$; $n=155,a=110,b=115$.
}
\label{simhyp_fig}
\end{figure}

\begin{figure}[!ht]
\begin{center}
\advance\leftskip-0.1cm
\includegraphics[width=4in]{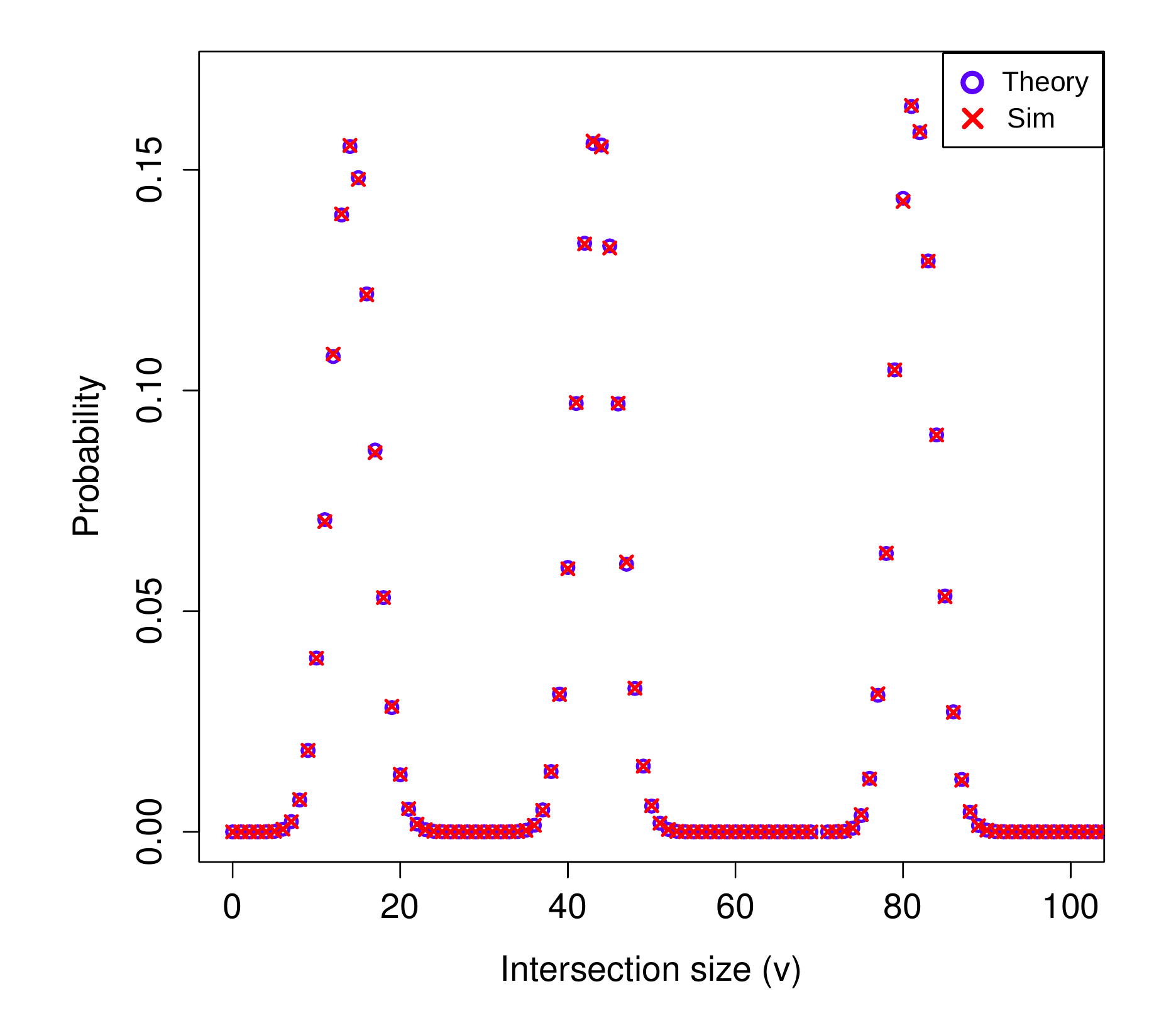}
\end{center}
\caption{
Match between theory and simulation for 3 parameter sets when sampling from 3 urns. Simulations (Sim) consisted of randomly and independently sampling without replacement from $n$ distinct categories in three separate urns and recording the size of the intersection each time (repeated 500,000 times for each distribution). From left to right, the parameters were: $n=100,a=57,b=41,c=61; n=100,a=83,b=76,c=69; n=140,a=123,b=119,c=109$.
}
\label{multi_3}
\end{figure}

\begin{figure}[!ht]
\begin{center}
\advance\leftskip-0.1cm
\includegraphics[width=4in]{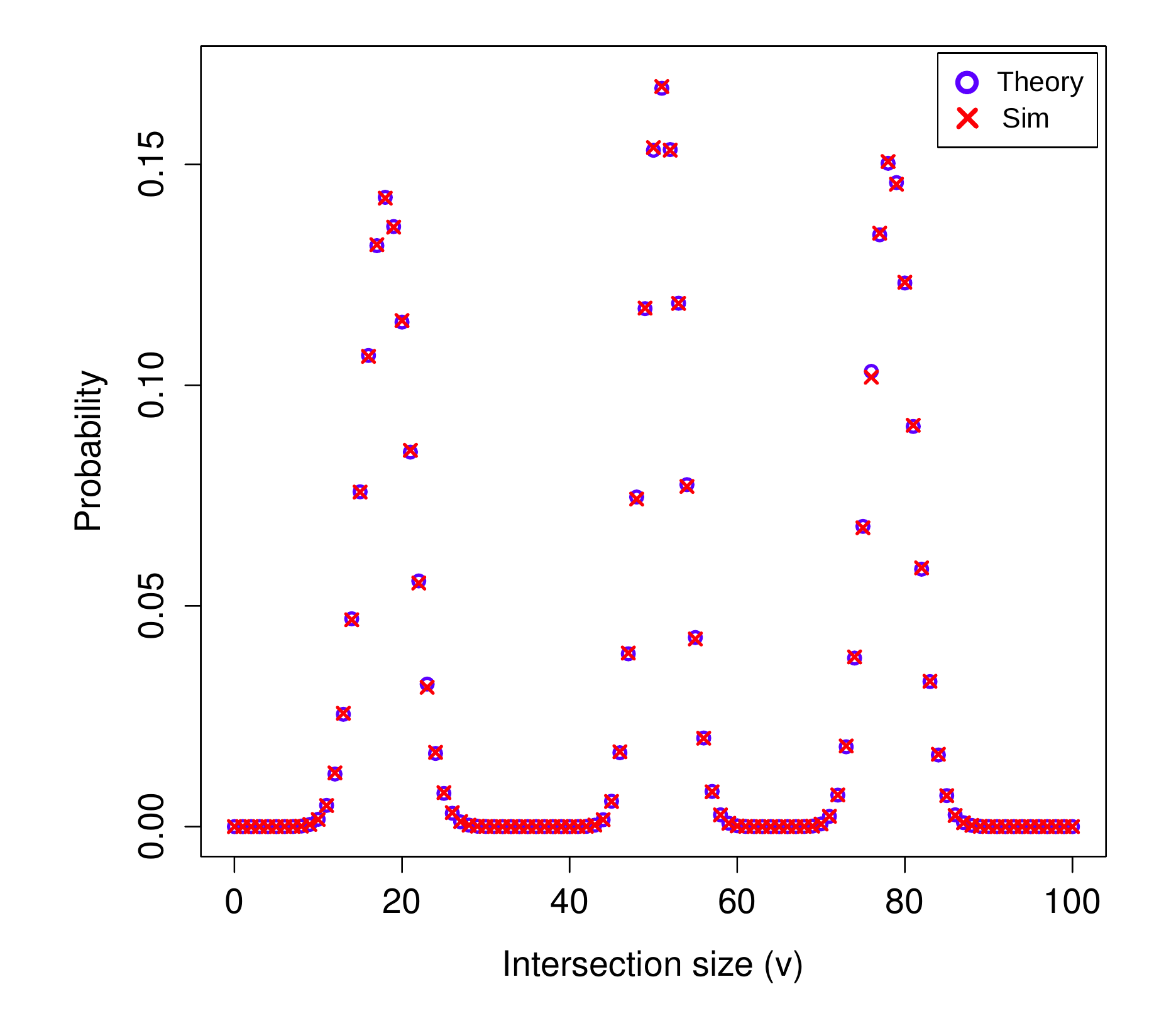}
\end{center}
\caption{
Match between theory and simulation for 3 parameter sets when sampling from 4 urns. Simulations (Sim) consisted of randomly and independently sampling without replacement from $n$ distinct categories in four separate urns and recording the size of the intersection each time (repeated 500,000 times for each distribution). From left to right, the parameters were: $n=100,a=64,b=79,c=58,d=62;n=100,a=92,b=81,c=77,d=89;n=140,a=121,b=118,c=131,d=115$.
}
\label{multi_4}
\end{figure}

\begin{figure}[!ht]
\begin{center}
\advance\leftskip-0.1cm
\includegraphics[width=4in]{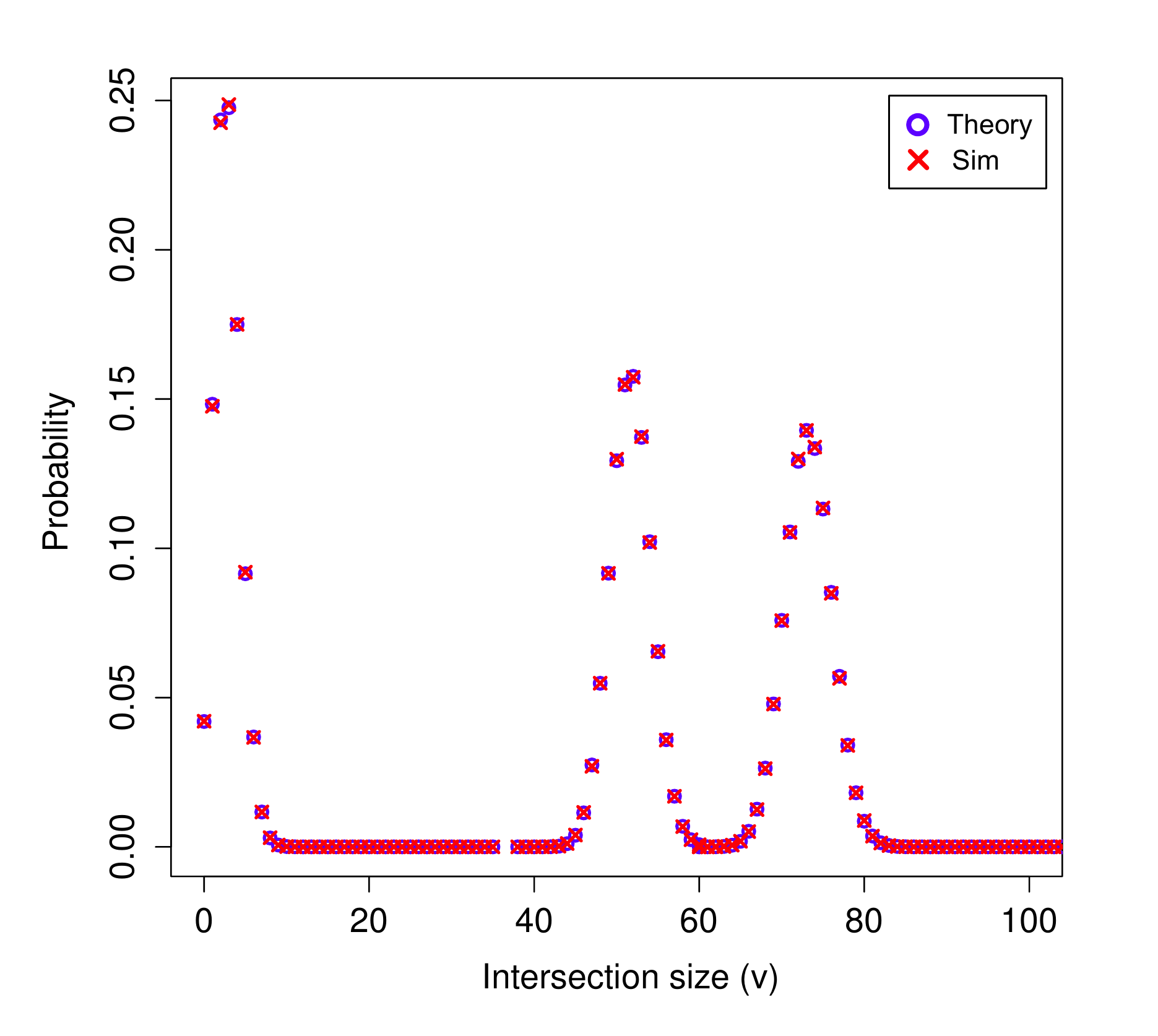}
\end{center}
\caption{
Match between theory and simulation for 3 parameter sets when sampling from 5 urns. Simulations (Sim) consisted of randomly and independently sampling without replacement from $n$ distinct categories in five separate urns and recording the size of the intersection each time (repeated 500,000 times for each distribution). From left to right, the parameters were: $n=108,a=35,b=43,c=84,d=63,e=49; n=101,a=85,b=93,c=84,d=91,e=89; n=138,a=122,b=118,c=119,d=126,e=123$.
}
\label{multi_5}
\end{figure}

\begin{figure}[!ht]
\begin{center}
\advance\leftskip-0.1cm
\includegraphics[width=4in]{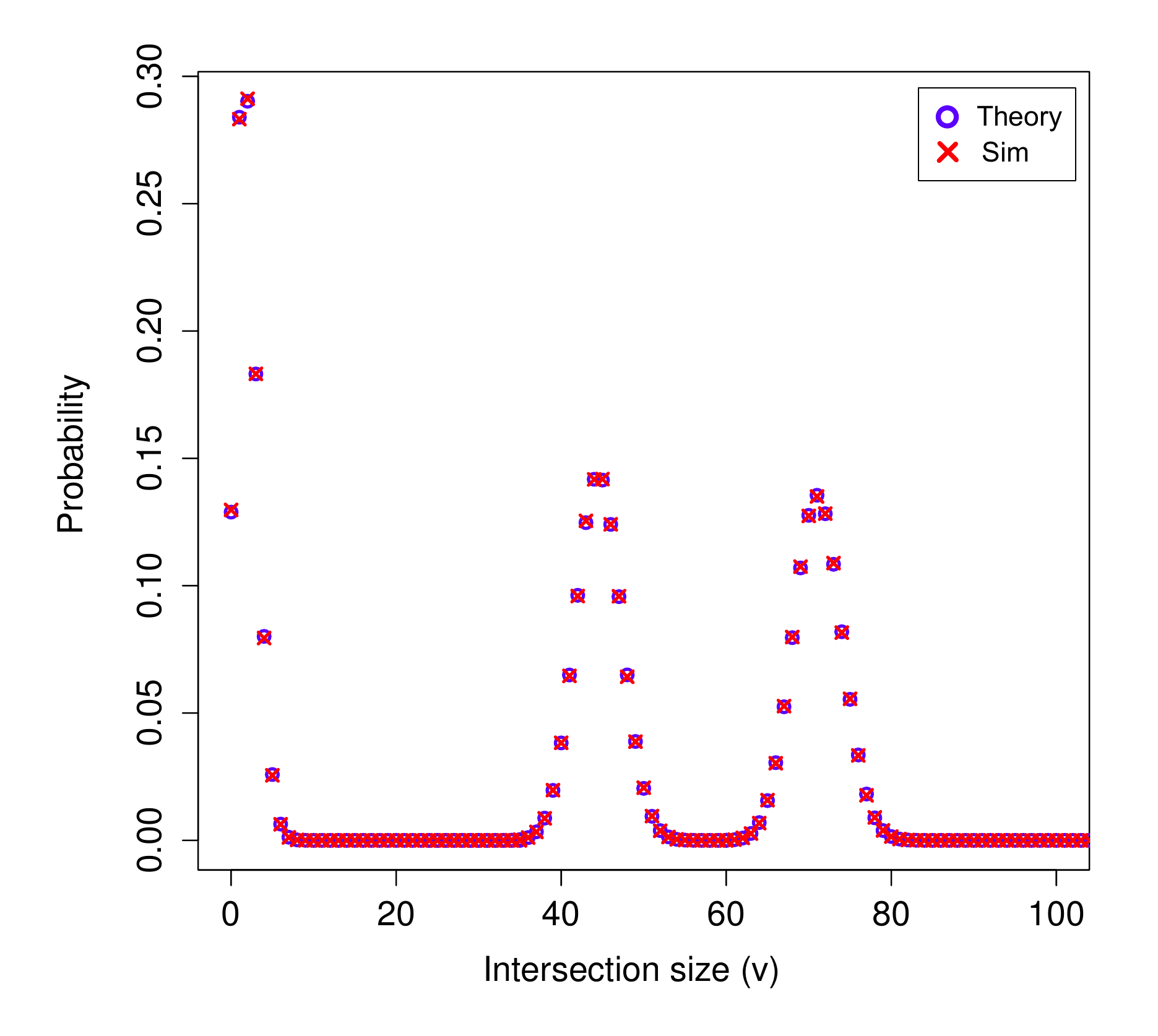}
\end{center}
\caption{
Match between theory and simulation for 3 parameter sets when sampling from 6 urns. Simulations (Sim) consisted of randomly and independently sampling without replacement from $n$ distinct categories in six separate urns and recording the size of the intersection each time (repeated 500,000 times for each distribution). From left to right, the parameters were: $n=n=108,a=35,b=43,c=84,d=63,e=49,f=72; n=101,a=85,b=93,c=84,d=91,e=89,f=87; n=138,a=122,b=118,c=119,d=126,e=123,f=134$.
}
\label{multi_6}
\end{figure}

\begin{figure}[!ht]
\begin{center}
\advance\leftskip-0.1cm
\includegraphics[width=4in]{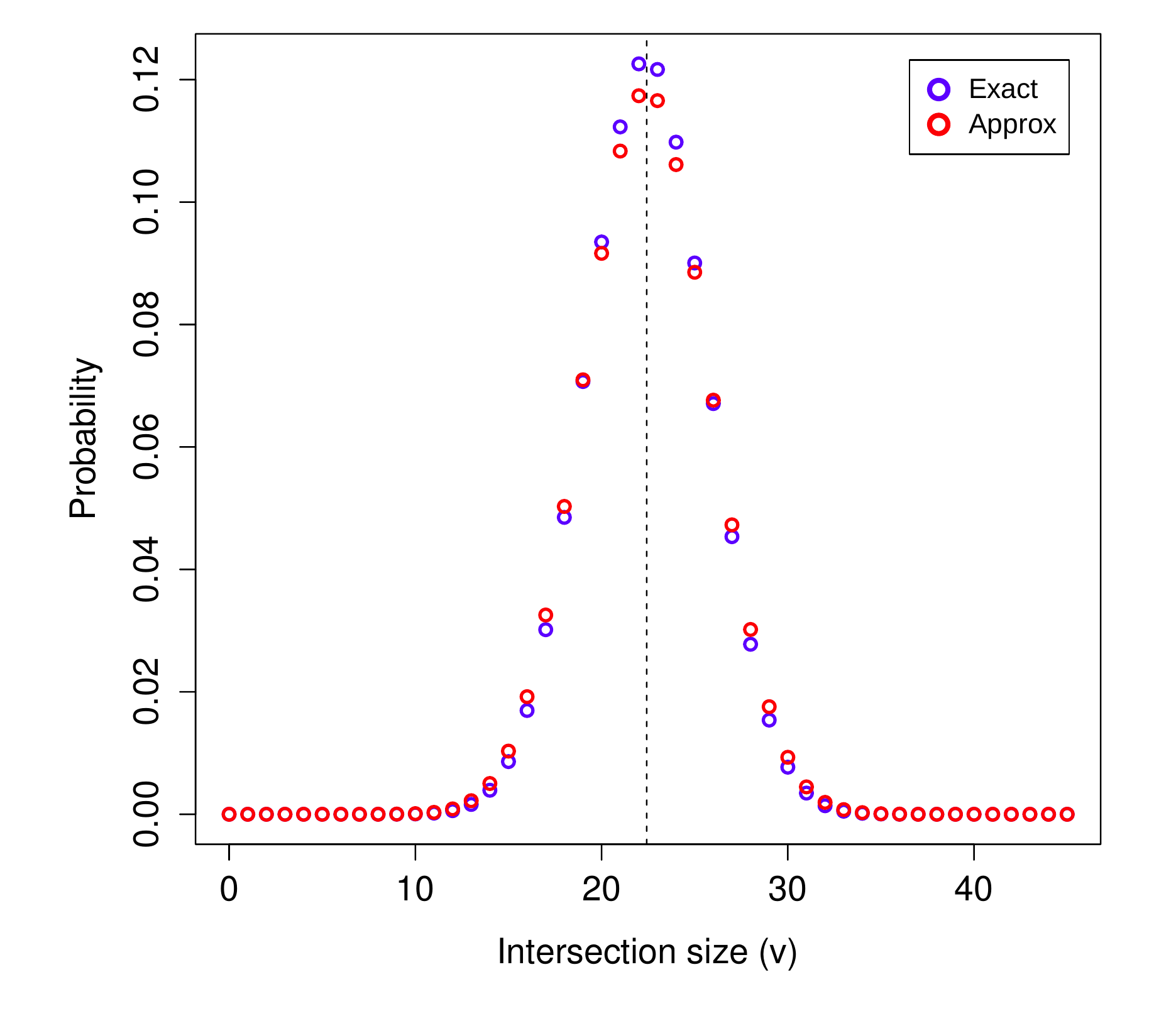}
\end{center}
\caption{
Match between true distribution and binomial approximation when sampling from 3 urns. The dashed vertical line indicates that the expectation is the same for both distributions. The parameters were: $n=452,a=361,b=45,c=282$.
}
\label{binom_approx}
\end{figure}

\begin{figure}[!ht]
\begin{center}
\advance\leftskip-0.1cm
\includegraphics[width=4in]{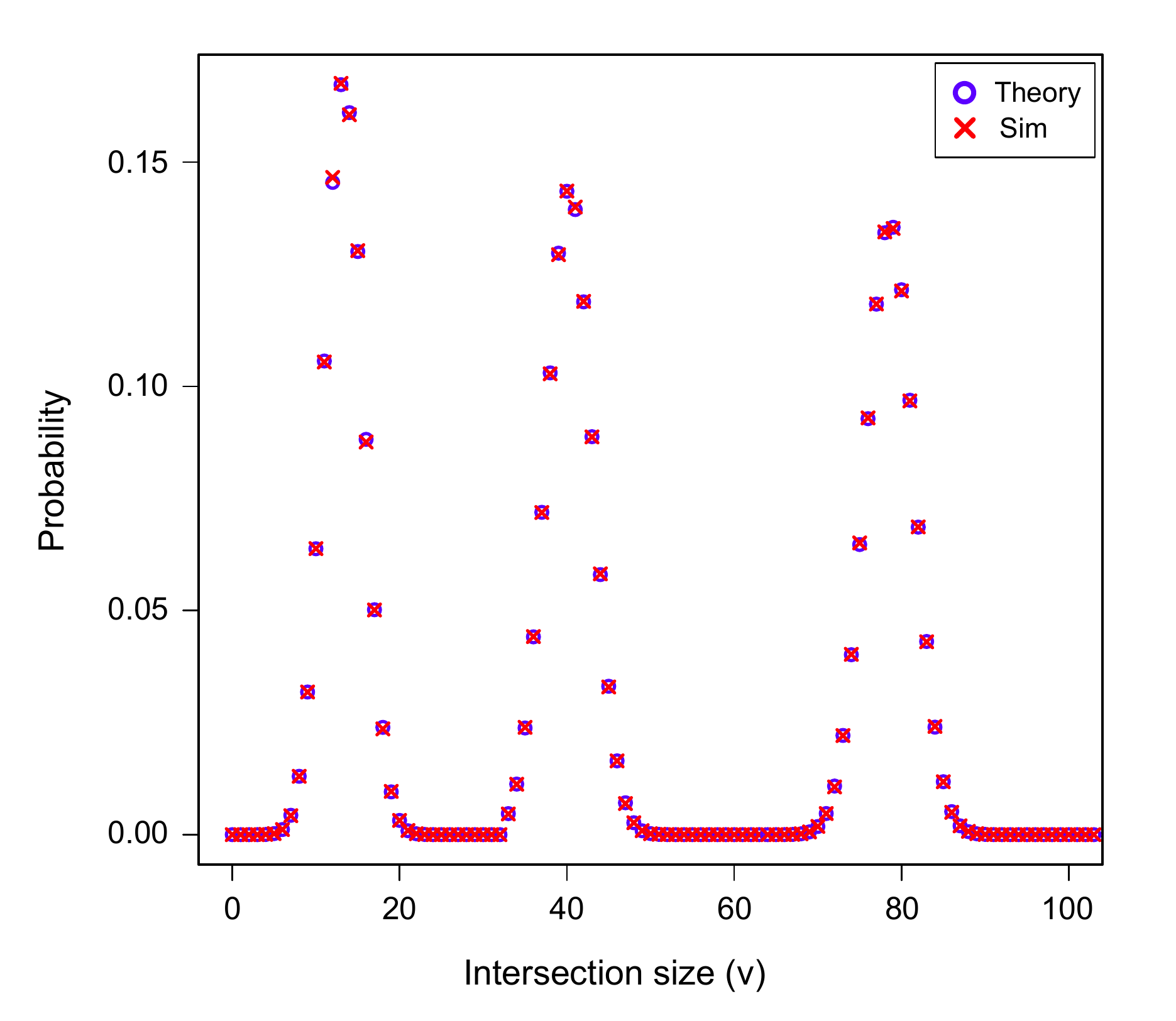}
\end{center}
\caption{
Match between theory and simulation for 3 parameter sets in the asymmetrical, duplicate case (\textbf{A} and \textbf{B$\mathbf{_2}$} in Figure \ref{urns_fig}). Simulations (Sim) consisted of randomly and independently sampling twice (without replacement) from $n$ distinct categories (in which the second set contained $q$ duplicates) and recording the size of the intersection each time (repeated 500,000 times for each distribution). From left to right, the parameters were: $n=100,a=35,b=42,q=59; n=100,a=63,b=79,q=73; n=130,a=110,b=115,q=47$.
}
\label{simhyp_dup_fig}
\end{figure}

\end{document}